\documentclass{scrartcl}

\usepackage{url} 

\usepackage{soul}

\usepackage{lmodern,lscape,textcomp,placeins}

\usepackage{lmodern,lscape,textcomp,placeins}
\usepackage[T1]{fontenc}

\usepackage{authblk}

\usepackage{amsmath,amssymb,dsfont,units,booktabs}
\usepackage{tabularx}
\usepackage{appendix}
\usepackage{graphicx, color}

\usepackage{graphicx}
\usepackage{lineno}
\usepackage{geometry}
\usepackage{xcolor}
\usepackage{setspace,amsmath}
\usepackage{bm}
\usepackage{authblk}
\usepackage[utf8]{inputenc} 
\newcommand{\bsigma}{\boldsymbol{\sigma}}
\newcommand{\Mop}{\mathsf{M}}
\newcommand{\Eop}{\mathsf{E}}
\newcommand{\Zop}{\mathsf{Z}}
\newcommand{\Sop}{\mathsf{S}}
\newcommand{\Top}{\mathsf{T}}
\newcommand{\Vop}{\mathsf{V}}
\newcommand{\Lop}{\mathsf{L}}

\newcommand{\Dop}{\mathsf{D}}

\date{\small\textbf{Abstract}}

\title{On discretizing sea-ice dynamics on triangular meshes using vertex, cell or edge velocities} 

{\small
\author[1,3,4]{S. Danilov}
\author[2]{ C. Mehlmann}
\author[1]{ V. Fofonova}
\affil[2]{Max-Planck Institute for Meteorology, Hamburg, Germany } 
\affil[1]{Alfred-Wegener-Institut, Helmholtz Zentrum für Polar- und Meeresforschung, Bremerhaven,
 Germany}
 \affil[3]{Jacobs University, Bremen, Germany} 
\affil[4]{A. M. Obukhov Institute of Atmospheric Physics RAS, Moscow, Russia}

}
\begin{document}

\maketitle

\begin{abstract}
 Discretization of the equations of Viscous Plastic and Elastic Viscous Plastic (EVP) sea ice dynamics on triangular meshes can be done by placing discrete velocities at vertices, cells or edges. Since there are more cells and edges than vertices, the cell- and edge-based discretizations simulate more linear kinematic features at the same mesh than the vertex discretization. However, the discretization based on cell and edge velocities suffer from kernels in the strain rate or stress divergence operators and need either special strain rate computations as proposed here for cell velocities, or stabilization as proposed earlier for edge velocities. An elementary Fourier analysis clarifies how kernels are removed, and also shows that cell and edge velocity placement leads to spurious branches of stress divergence operator with large negative eigenvalues. Although spurious branches correspond to fast decay and are not expected to distort sea ice dynamics, they demand either smaller internal time steps or higher stability parameters in explicit EVP-like methods.   
\end{abstract}

\section{Introduction}
Sea-ice, located at high-latitudes and at the boundary between ocean and atmosphere, plays an important role in the climate system. 
Modelling the complex mechanical and thermodynamical behaviour of sea-ice at a broad range of spatio-temporal scales poses a manifold of challenges. 
Freezing sea water forms a composite of pure ice, liquid brine, air pockets and solid salt. The details of this formation depend on the laminar or turbulent environmental conditions. 

\section{Introduction}
Currently  almost all sea ice models treat sea ice as a viscous-plastic material  either in the framework of the viscous plastic (VP) rheology of \cite{Hibler1979} or by using the elastic viscous plastic (EVP) formulation introduced by \cite{HunkeDukowicz1997}.

Even though the use of VP rheology at a grid spacing of the size of a single sea ice floe is questioned (e.g. \cite{Coon2007,Feltham08}), recent work of \cite{Wang2016} and \cite{Hutter2018a} indicates that sea ice models based on VP rheology simulate many highly localized deformation features that are observed by Synthetic Aperture Radars (SAR). These features, which  are referred to as linear kinematic features (LKFs), start to appear in simulation with a gird spacing of about 4 km (\cite{Wang2016}). While the effect of simulated LKFs on the exchange between the ocean and atmosphere is an emerging topic for further research, a question also arises on how the simulated LKFs are related to discretization of sea ice dynamical equations. The study of \cite{Mehlmannetal2021} compares several discretizations on quadrilateral and triangular meshes, showing that the placement of the sea ice velocity plays a major role in defining the number of simulated LKFs for a given mesh resolution.

The discretization of sea-ice dynamics was traditionally done on quadrilateral meshes assuming Arakawa B- or C-grid placement (\cite{Arakawa1977}) for sea ice velocities and scalar variables. The appearance of global ocean circulation models formulated on unstructured meshes such as FESOM (\cite{Wang2014}, \cite{Danilov2017}), MPAS-O (\cite{Ringler2013}) and ICON-O (\cite{Korn2017}), as well as the need for a sea ice component in coastal ocean models (\cite{Gao2011}) raised a question on reformulating the discrete equations of sea ice dynamics on unstructured (triangular or dual hexagonal) meshes. This task was also considered by \cite{Hutchings2004} for general polygonal meshes and \cite{Lietaer2008} for triangular meshes outside the framework of a specific ocean model. 

The three recent basic approaches to discretize sea ice dynamics on triangular and hexagonal meshes differ by the placement of discrete degrees of freedom. The sea-ice component of FESOM (\cite{Timmermann2009,Danilov2015}) relies on vertex placement of ice velocities and scalars (e.g. sea ice concentration and thickness), and was implemented using $P_1$ linear continuous finite elements. It corresponds to A-grid discretization, and can be also seen as a finite-volume discretization for median-dual control volumes around mesh vertices. \cite{Hutchings2004} place horizontal velocities and scalars on mesh cells. \cite{Gao2011} also take sea ice velocity vectors on cells, but use vertex-based scalars, which is an analog of B-grid (because of staggering and full velocity vectors). The MPAS-O discretization as described in \cite{Petersen2019} is analogous to triangular B-grid. Finally \cite{Lietaer2008} and recently \cite{Mehlmann2021} propose to use the edge placement of sea ice velocity vectors, treated with linear non-conforming (Crouzeix--Raviart (CR)) finite elements. The scalars are  constant on cells in both cases. This placement corresponds to a CD-grid type staggering.

On large triangular meshes the numbers of vertices, cells and edges are related as 1:2:3. This is why setups with discrete velocities placed on triangles or edges have more velocity degrees of freedom (DoF) than setups based on vertex velocities. They may ensure a better spatial resolution on the same mesh. Indeed,  \cite{Mehlmannetal2021}  demonstrated that on triangular meshes the number of simulated LKFs is the highest for discretizations placing velocity vectors on edges (CD-grid), followed by the placement on cell centers (B-grid). They both outperform the vertex-based collocated discretization (A-grid). 

However, in contrast to the A-grid discretization, the cell and edge placements support numerical modes related to the geometry of triangular mesh. Furthermore, strain rates or stress divergence for the cell and edge placements may possess kernels, i.e. be zero for non-trivial discrete velocities. The intention of present work is to clarify the origin of numerical modes and kernels, and discuss measures allowing to handle them. They will be referred to as stabilization. They have already been used in FESOM in the benchmark comparison of \cite{Mehlmannetal2021}, but without theoretical analysis. 

The instability of a CD-grid discretization of the sea ice dynamics on triangular meshes was recently mentioned by \cite{Mehlmann2021} who proposed a stabilization that controls oscillations by removing the kernel in discrete stress divergence. Although the instabilities accompanying the cell placement of velocity on triangular meshes were not documented in the literature cited above, the discretizations proposed in \cite{Hutchings2004} and \cite{Gao2011} as well as a triangular-mesh analog of the discretization in \cite{Petersen2019} have kernels either in discrete strain rates or stress divergence. We explain how kernels are created in this case and propose a stable B-grid discretization. We also show that numerical modes supported on B- and CD-grids correspond to anomalously large negative eigenvalues of stress divergence operator, which has implications for explicit time stepping in the EVP method or its modified version mEVP (\cite{Bouillon2013}).
 
We use a standard Fourier analysis in section \ref{sec:fourier} as a basic tool to analyze the behavior of A-, B- and CD-grid discretizations. Our analysis is therefore limited to regular triangular meshes and linear (viscous) regimes. Despite these limitations, we hope that it adds to the understanding of the core difficulties. 
A test case proposed in \cite{Mehlmannetal2021} is used in section \ref{sec:ill} to illustrate the absence of spurious oscillations in the stabilized setups and to analyze the time step restriction of different discretizations. The equations of VP, EVP and mEVP dynamics can be found in papers cited above (e.g., \cite{Danilov2015}).  

\section{Fourier analysis of stress divergence} \label{sec:fourier}
\subsection{Velocity representation}
The analysis focuses on the vertex, cell and edge  placements for discrete sea ice velocities on triangular meshes. The three options correspond to the velocity staggering of the A-grid, B-grid and CD-grid discretizations respectively. 

We consider mesh made of equilateral triangles and introduce a coordinate system with axes $x$ and $y$. We may select two triangles, one with vertices at $(0,0)$, $(a,0)$ and $(a/2,h)$, and the other one with vertices at $(a,0)$, $(3a/2,h)$ and $(a/2,h)$. Here $a$ is the side of triangle and $h=\sqrt{3}a/2$ is the height. All mesh triangles are obtained by translations of these two along their sides, i.e., by the displacements 
$$
\mathbf{q}=n_1(1,0)a+n_2(1/2,\sqrt{3}/2)a,
$$
where $n_1$ and $n_2$ are integers. We will refer to triangles of the first type as $u$ (pointing up) triangles, and to triangles of the second type as $d$ (pointing down) triangles. The mesh is the union of two sets $\mathcal{C}^u$ and $\mathcal{C}^d$ of triangles. 

Although sea ice stresses $\bsigma$ of VP rheology (\cite{Hibler1979}) depend nonlinearly on velocities, this dependence becomes linear in the viscous regime. This is the only case which allows a Fourier analysis. We assume for the rest of this section that we are in the viscous regime, i.e. we deal with constant viscosities $\eta$ and $\zeta$ in the expression for the components of $\bsigma$: 
$$
\sigma_{ij}=2\eta(\dot\epsilon_{ij}-\frac{1}{2}\delta_{ij}\dot\epsilon_{nn})+\zeta\delta_{ij}\dot\epsilon_{nn}-\delta_{ij}\frac{P}{2},\,\,\dot\epsilon_{ij}=\frac{1}{2}(\partial_i u_j+\partial_ju_i),\, i,j=x,y,
$$
where $\mathbf{u}=(u_x,u_y)=(u,v)$ is the sea ice velocity, $P$ is the ice strength, $\delta_{ij}$ is the Kronecker delta and summation is implied over repeating coordinate indices. In the VP rheology, $\eta=\zeta/e_{VP}^2$, with $e_{VP}$ the ratio of major to minor axes of the elliptical yield curve, and $\zeta=P/(2\max(\Delta,\Delta_{min}))$. Here $\Delta=((\dot\epsilon_{xx}+\dot\epsilon_{yy})^2+(1/e^2_{VP})((\dot\epsilon_{xx}-\dot\epsilon_{yy})^2+4\dot\epsilon_{xy}^2))^{1/2}$, and $\Delta_{min}$ is the parameter defining transition between plastic and viscous regimes. The viscous regime takes place if $\Delta>\Delta_{min}$. We will assume $P=\mathrm{const}$ to ensure  the constancy of $\eta,\zeta$. In this case the $P$ term will not contribute to the stress divergence. 

We will be interested in the behavior of discrete divergence of stresses: \begin{equation}
\Vop\mathbf{u}=\nabla\cdot\bsigma.
\label{eq:sigdiv}    
\end{equation} 
In the continuous case, taking $\mathbf{u}=(\overline{u}, \overline{v}) e^{i\mathbf{k}\cdot\mathbf{x}}$, where $\overline{u}$ and $\overline{v}$ are the amplitudes of $x$ and $y$ components of velocity, and inserting in (\ref{eq:sigdiv}), we see that the Fourier symbol of $\Vop$ is a 2 by 2 matrix which returns the amplitudes of $\Vop\mathbf{u}$ if applied to the vector of velocity amplitudes $(\overline{u}, \overline{v})^T$. The eigenvalues of this matrix are $-\eta(k^2+l^2)$ and $-(\eta+\zeta)(k^  
2+l^2)$, where $k$ and $l$ are the $x$- and $y$-components of the wavevector $\mathbf{k}$. The eigenvalues correspond, respectively, to the transverse and longitudinal (with respect to $\mathbf{k}$) modes. We will discuss further the extent to which these eigenvalues are modeled by discrete solutions.  Although in the case of VP rheology the ratio $z=\zeta/\eta=e^2_{VP}$, with $e_{VP}=2$, we will use $z=1$ in this section to make the eigenvalues closer to each other in graphical representation.

The form of the Fourier solution in the discrete case depends on the placement of velocities. Since mesh vertices are invariant to the set of translations $\mathbf{q}$, the discrete vertex velocities will be taken in the form
$$
\mathbf{u}_v=\overline{\mathbf{u}}e^{i\mathbf{k}\cdot\mathbf{x}_v},\quad v\in\mathcal{V},
$$     
where $\overline{\mathbf{u}}$ is the vector of amplitudes, $\mathbf{x}_v$ is the vector drawn to vertex $v$, and $\mathcal{V}$ is the set of mesh vertices. In the case of cell velocities the set of translations $\mathbf{q}$ connects only subsets of $u$ and $d$ triangles. We need, therefore, separate vectors of amplitudes for $u$ and $d$ triangles, and we search for discrete solutions in the form 
$$
\mathbf{u}^u_c=(\overline{u}^u,\overline{v}^u)e^{i\mathbf{k}\cdot\mathbf{x}_c},\quad c\in\mathcal{C}^u,
$$    
$$
\mathbf{u}^d_c=(\overline{u}^d,\overline{v}^d)e^{i\mathbf{k}\cdot\mathbf{x}_c},\quad c\in\mathcal{C}^d.
$$
Here $\mathbf{x}_c$ is the vector drawn to the center of triangle $c$.
Finally, three vectors of amplitudes will be needed in the case of edge velocities   
$$
\mathbf{u}^a_e=(\overline{u}^a,\overline{v}^a)e^{i\mathbf{k}\cdot\mathbf{x}_e},\quad e\in\mathcal{E}^a,
$$    
$$
\mathbf{u}^b_e=(\overline{u}^b,\overline{v}^b)e^{i\mathbf{k}\cdot\mathbf{x}_e},\quad e\in\mathcal{E}^b,
$$
and
$$
\mathbf{u}^c_e=(\overline{u}^c,\overline{v}^c)e^{i\mathbf{k}\cdot\mathbf{x}_e},\quad e\in\mathcal{E}^c,
$$
where $\mathbf{x}_e$ is the vector drawn to the center of edge $e$. Here the set $\mathcal{E}$ of mesh edges is split into the three subsets $\mathcal{E}^a,\mathcal{E}^b$ and $\mathcal{E}^c$ of edges oriented as $(a,0)$, $(a/2,\sqrt{3}a/2)$ and $(-a/2,\sqrt{3}a/2)$ respectively. The Fourier symbol of $\Vop\mathbf{u}$ will be a $2\times2$ matrix in the case of vertex velocities, because there are only two Fourier amplitudes, as in the continuous case. It will be a $4\times4$ matrix in the case of cell placement, acting on the vector $(\overline{u}^u,\overline{v}^u,\overline{u}^d,\overline{v}^d)^T$, and a $6\times6$ matrix in the case of edge velocities, acting on the vector $(\overline{u}^a,\overline{v}^a,\overline{u}^b,\overline{v}^b,\overline{u}^c,\overline{v}^c)^T$. There will be, respectively, four and six branches in the cell and edge cases instead of two. Commonly the presence of extra branches implies that some of them will be spurious, and the analysis below indicates that it is indeed so for cell and edge velocities. 

\subsection{Vertex velocities}
\label{sec:p1}
\begin{figure}[ht]
\centering
  \includegraphics[width=12cm]{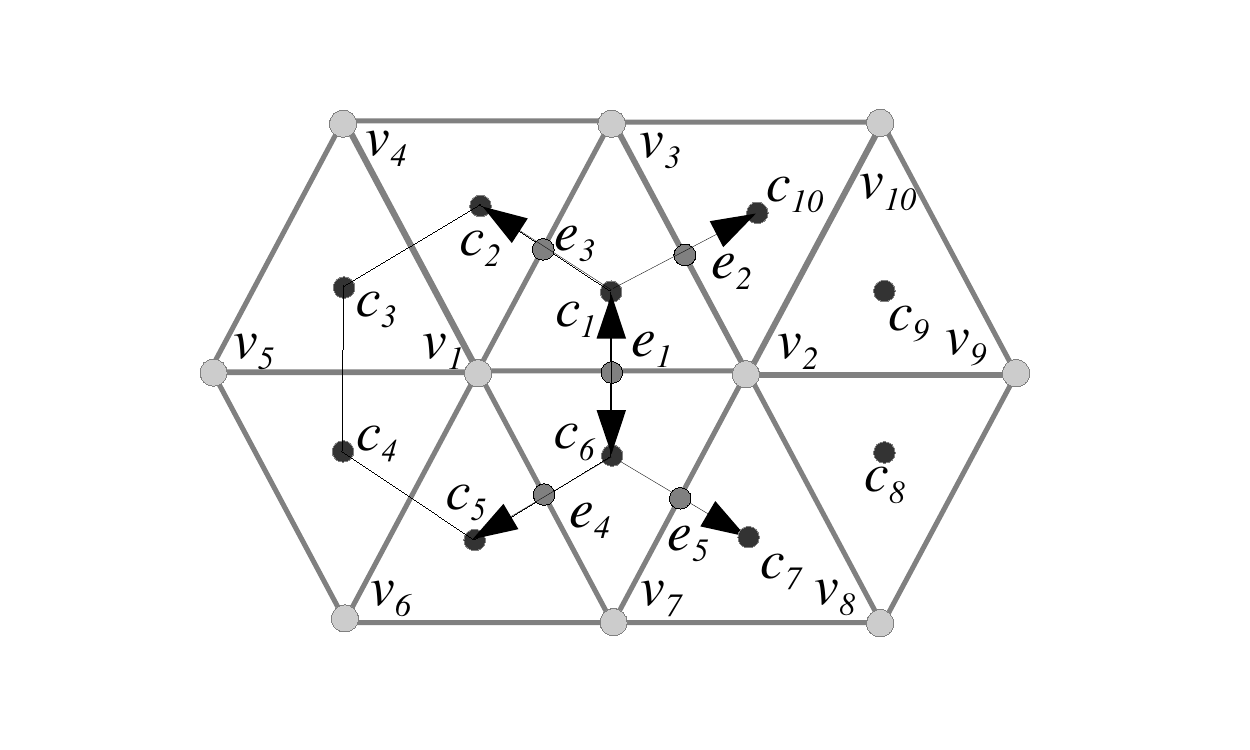}
\vspace{-0.5cm}\caption{Schematic of mesh geometry. For $P_1$ vertex velocities, strain rates and stresses are computed at triangles (e.g., $\mathbf{u}_{v_1}$, $\mathbf{u}_{v_2}$ and $\mathbf{u}_{v_3}$ define the strain rates at $c_1$) and stress divergence returns the result to vertices according to (\ref{eq:p1a}). For cell velocities, in case V, strain rates and stresses are computed at vertices using Gauss' theorem for median-dual control volumes around vertices. A median-dual control volume shown for vertex $v_1$ is formed by joining centers of triangles containing $v_1$ with centers of edges containing $v_1$ (thin lines). The divergence of stresses is computed on triangles assuming that the stresses are linear on triangles. In case C strain rates are computed at cells using least squares linear fit over neighbors (e.g. $c_1,c_2,c_6$ and $c_{10}$ for $c_1$). Stresses are computed on cells and then averaged to edges (e.g. the stresses at $e_1$ are the half sum of stresses at $c_1$ and $c_6$). Stress divergence is computed by applying Gauss' theorem to triangles. For edge velocities strain rates are at triangles, and the divergence of stresses is at edges (e.g., strain rates at $c_1$ get contributions from $e_1, e_2$ and $e_3$, and the strain divergence at $e_1$ gets contribution from $c_1$ and $c_6$ in the absence of stabilization). }
\label{fig:schematic}
\end{figure}

We use finite-element discretization (see \cite{Danilov2015}) and compute the discrete $(\Vop \mathbf{u})_v$ writing
\begin{equation}
\int\mathbf{N}_v\cdot(\Vop \mathbf{u}) dS=-\int\nabla \mathbf{N}_v:\bm{\sigma} dS,
\label{eq:p1a}
\end{equation}
where $\mathbf{N}_v=\mathbf{w}_vN_v$ is the vector test function, with the amplitude equaled to amplitude of $\mathbf{w}_v$ at vertex $v$, and $N_v$ the standard $P_1$ linear function equal to 1 at vertex $v$, decaying linearly to 0 at neighbor vertices and being 0 outside the stencil of triangles containing $v$. 
The discrete sea ice velocities and stress divergence are expressed as
$$
\mathbf{u}=\sum_{v'\in\mathcal{V}} N_{v'}\mathbf{u}_{v'}, \quad
\Vop\mathbf{u}=\sum_{v'\in\mathcal{V}} N_{v'}(\Vop\mathbf{u})_{v'}.
$$
If the mass matrix $\int N_vN_{v'}dS$ appearing on the left hand side of (\ref{eq:p1a}) is lumped, its diagonal $v$-entry $A_v$ is the area of median-dual control volume around vertex $v$, i.e.,  $A_v=\sum_{c\in\mathcal{C}(v)} A_c/3$. Here, $\mathcal{C}(v)$ is the set of triangles containing $v$ and $A_c$ is the area of triangle $c$. $A_v$ is twice the triangle area if the mesh is uniform. For mesh patch in Fig. \ref{fig:schematic}, the stress divergence $(\Vop \mathbf{u})_v$ at $v=v_1$ will be defined by contributions from triangles $c_1-c_6$ if the mass matrix is lumped.

For cell-wise linear velocities, velocity derivatives (and hence strain rates and stresses) are constant on triangles. Although these constant values are related to the entire triangles, we will interpret them as located at centers of triangles, which is second-order accurate. This is done only to have a rule to compute the phase multiplier in the expression like
$$
(\nabla u)_c=(\overline{\nabla u})e^{i\mathbf{k}\cdot\mathbf{x}_c}=\sum_{v\in\mathcal{V}(c)} u_v\nabla N_v=\left(\sum_{v\in\mathcal{V}(c)} \overline{u}e^{i\mathbf{k}\cdot(\mathbf{x}_v-\mathbf{x}_c)}\nabla N_v\right)e^{i\mathbf{k}\cdot\mathbf{x}_c}=
(\mathbf{G}^u\overline{u})e^{i\mathbf{k}\cdot\mathbf{x}_c}
$$
for $c\in\mathcal{C}^u$, and likewise for the $v$ component.
Here, $\mathcal{V}(c)$ is the set of vertices of triangle $c$, and
$$
\mathbf{G}^u=(g_x,g_y)=-\frac{1}{2h}\left( (\sqrt{3},1)\alpha_1+(-\sqrt{3},1)\beta_1+(0,-2)\gamma_1\right),
$$
where $\quad\alpha_1=e^{-ika/2-ilh/3},\, \beta_1=e^{ika/2-ilh/3},\,\gamma_1=e^{2ilh/3}$.
For $c\in\mathcal{C}^d$, we will deal with $\mathbf{G}^d=-(g_x^*, g_y^*)$, the star implies complex conjugation. It can be readily seen that the Fourier symbols $g_x$ and $g_y$ only include  phase differences between the vertices of triangle and its center, and that these differences depend on the kind ($u$ or $d$) of triangle. As the consequence, we can compute matrices of Fourier symbols for strain rates and stresses on $u$ and $d$ triangles. The strain rate amplitudes are connected to velocity amplitudes as 
$$
(\overline{\dot\epsilon_{xx}^u},\overline{\dot\epsilon_{xy}^u},\overline{\dot\epsilon_{yy}^u},\overline{\dot\epsilon_{xx}^d},\overline{\dot\epsilon_{xy}^d},\overline{\dot\epsilon_{yy}^d})^T=\Eop\begin{pmatrix}\overline{u}\\ \overline{v}\end{pmatrix}=\begin{pmatrix}g_x&0\\g_y/2 &g_x/2\\0 &g_y\\ -g_x^*&0\\-g_y^*/2 &-g_x^*/2\\0 &-g_y^*\end{pmatrix}\begin{pmatrix}\overline{u}\\ \overline{v}\end{pmatrix}.   
$$
The expression connecting the strain rates and stresses corresponds to the matrix     
\begin{equation}
    \label{eq:Smat}
\Zop=\begin{pmatrix}\Sop & \mathbf{0}\\\mathbf{0}& \Sop\end{pmatrix},\quad
\Sop=\eta\begin{pmatrix}1+z & 0& -1+z\\0&2&0\\-1+z&0&1+z\end{pmatrix}.
\end{equation} 
Finally, the divergence of stresses on the right hand side of (\ref{eq:p1a}), is the sum of contributions from six triangles. One readily sees that the contributions from $u$ triangles come with opposite sign and phases compared to $g_x$ and $g_y$ (the centers of $u$-triangles around vertex $v$ form a $d$ triangle), giving     
$$
\Dop=\frac{1}{2}\begin{pmatrix}-g_x^*&-g_y^*&0&g_x&g_y&0\\ 
 0&-g_x^*&-g_y^*&0&g_x&g_y\end{pmatrix}.
$$
$\Dop$ already incorporates the division by $A_v$ so that $\Dop\Zop\Eop$ is the Fourier symbol of $\Vop$ if the mass matrix is lumped. The Fourier symbol of consistent mass matrix is 
$$
A_v\Mop=A_v\begin{pmatrix} m&0\\0&m\end{pmatrix},\quad m=1/2+(\cos(ka)+\cos(ka/2+lh)+\cos(-ka/2+lh))/6,
$$
leading to the expression $\Mop^{-1}\Dop\Zop\Eop$ for the Fourier symbol of $\Vop$.

\begin{figure}[t]
\centering
  \includegraphics[width=12cm]{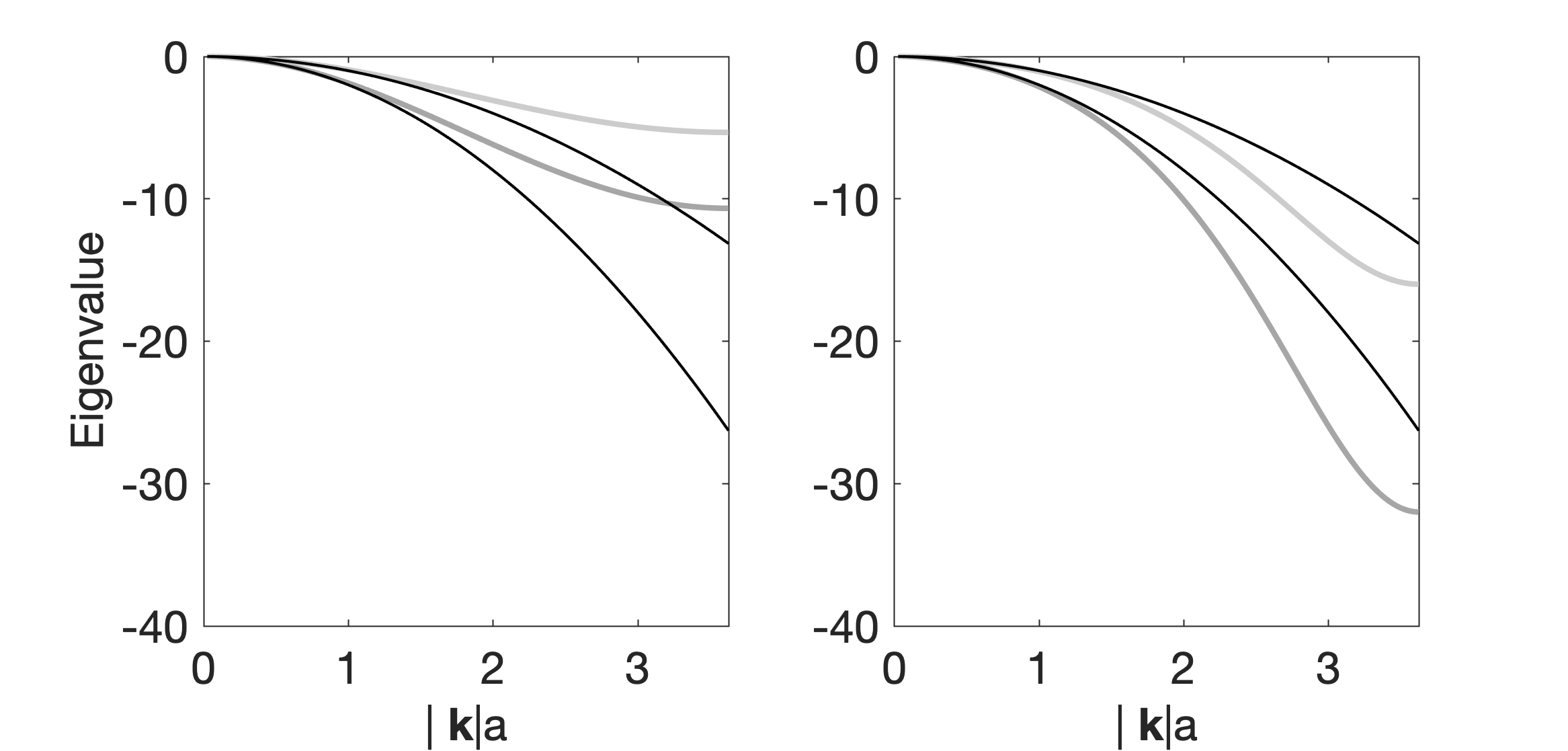}
\vspace{-0.5cm}\caption{The eigenvalues of the Fourier symbol of $a^2\eta^{-1}\nabla\cdot\bsigma$ (thick gray lines) for vertex velocity placement as a function of dimensionless wavenumber $|\mathbf{k}|a$ for the wavevector oriented at $\pi/6$ to the $x$-axis and $\zeta=\eta$. The thin black lines plot the dimensionless eigenvalues $-(k^2+l^2)a^2$ and $-2(k^2+l^2)a^2$ of the continuous case. Left panel: lumped mass matrix, Right panel: consistent mass matrix.} 
\label{fig:p1}
\end{figure}

The left panel of Fig. \ref{fig:p1} presents the eigenvalues of $a^2\eta^{-1}\Dop\Zop\Eop$ as a function of $|\mathbf{k}|a$ for $\mathbf{k}$ directed at $\pi/6$ to the $x$-axis. The boundary of the first Brillouin zone for a triangular mesh composed of equilateral triangles in this direction is at $|\mathbf{k}|a=2\pi/\sqrt{3}$ (see, e.g., \cite{DanilovKutsenko2019}). It defines the largest resolvable wavenumber. The eigenvalues are close to the theoretically predicted ones if wavenumbers are sufficiently small. This behavior is preserved for all directions of the wavevector, with some spread of curves only in the region of largest wavenumbers. One readily sees that the operator is sufficiently accurate at scales larger than $3h$. From the right panel, which shows the eigenvalues for the consistent mass matrix, one may conclude that the accuracy is not improved if a consistent mass matrix is used, the error only changes the sign. This means that using a lumped mass matrix is fully appropriate for vertex velocities.    

\subsection{Cell velocities}

First described by \cite{Hutchings2004} for arbitrary polygonal meshes, the cell placement of sea ice velocities on triangular meshes is used in FVCOM (\cite{Gao2011}). Its analog on hexagonal meshes is the B-grid discretization in MPAS-O (\cite{Petersen2019}), where the circumcenters of triangular cells are the corners of hexagons of dual mesh. This placement was also explored for FESOM2 setup (\cite{Danilov2017}), but found to be leading to noise and was abandoned in favor of older vertex placement. \cite{Hutchings2004} and \cite{Gao2011} compute strain rates on cells centers using a least squares fit on the stencil of nearest cells, with subsequent averaging to edges.  \cite{Petersen2019} computed the strain rates at centers of hexagonal cells, which correspond to vertices on triangular meshes, using a variational approach. The cell and vertex computation of strain rates will be referred to as C and V cases (see Fig. \ref{fig:schematic}). Unfortunately, both face instabilities on triangular (or dual) meshes related to the mesh geometry, and need adjustments. 

\subsubsection{Strain rate computation}
We begin with the V case. Using the generalized form of Gauss' theorem, the velocity derivative $\partial_x u$ can be computed as
$$
A_v(\partial_x u)_v=\sum_{s\in \mathcal{S}(v)}\ell_s u_s (n_x)_s,
$$
and similarly for the remaining components. In this expression, $s$ denotes segments of the boundary of median-dual control volume around vertex $v$ (or the boundary of hexagonal cell with center $v$) as shown schematically in Fig. \ref{fig:schematic}, $\mathcal{S}(v)$ is the set of such segments for particular $v$, $\ell_s$ is the length of segment, $(n_x)_s$ is the $x$-component of the segment outer normal.  On general meshes, the median-dual control volumes will differ from the cells of dual mesh, but they coincide on the regular equilateral mesh. Since the boundary of the control volume around $v$ passes through the cell centers, the attribution of cell velocities to segments is trivial.

The velocity derivatives are used to compute the components $\dot\epsilon_{ij}$ of symmetric strain rate tensor. If a mesh is made of $N$ vertices, there are approximately $2N$ triangles if the mesh is large enough to neglect the effect of boundaries, and hence $4N$ degrees of freedom in sea ice velocities. However, strain rates and stresses on vertices have only $3N$ degrees of freedom. The rank of the $4N\times 4N$ matrix expressing the discrete divergence of stresses in terms of velocities is not larger than $3N$, i.e., a nontrivial null-space is created. Since the null-space is created by computing strain rates, it will persist beyond the linear viscous regime.  

In case C (see schematic in Fig.\ref{fig:schematic}), stresses (and strain rates) are computed at triangles and then averaged to edges before computing the stress divergence. On uniform meshes, the contributions from the nearest triangles drop out from the expression for the stress divergence in this procedure, allowing a mode in velocities. In Fig. \ref{fig:schematic}, stress divergence at $c_1$ will miss the contributions from velocities at $c_2,c_6$ and $c_{10}$. The implicit computation of gradients at edges in \cite{Hutchings2004} hints at the method that can be used to eliminate the mode. 

We now illustrate this with the Fourier analysis, explaining also an approach that is free of difficulties.

\subsubsection{Fourier analysis for cell velocities}
For the analysis we associate the cell velocities with cell centers because V and C cases assume this to compute stresses. In case V the same geometry of cells and vertices is involved as in the computations of stress divergence for the vertex velocities in section \ref{sec:p1}. The strain rates will be defined by the velocities at three $u$ and three $d$ triangles meeting at $v$ ($c_1, c_3, c_5$ and $c_2, c_4, c_6$ for $v_1$ in Fig. \ref{fig:schematic}). We will get 
$$
\overline{\dot\epsilon_{xx}}=(1/2)(-g_x^*,0,g_x,0)(\overline{u}^u,\overline{v}^u,\overline{u}^d,\overline{v}^d)^T  
$$
for the amplitude of $\dot\epsilon_{xx}$, and similarly for the other strain rate components. 
The matrix
$$
\Eop=\frac12\begin{pmatrix}-g_x^* &0&g_x&0\\ -g_y^*/2&-g_x^*/2&g_y/2&g_x/2\\
0&-g_y*&0&g_y\end{pmatrix}
$$ 
is the Fourier symbol of strain rates connecting the vector of their amplitudes $(\overline{\dot\epsilon_{xx}}, \overline{\dot\epsilon_{xy}}, \overline{\dot\epsilon_{yy}})^T$ with the velocity amplitudes $(\overline{u}^u, \overline{v}^u, \overline{u}^d, \overline{v}^d)^T$. The amplitude of stresses are connected to the amplitudes of strain rates by the matrix $\Sop$ in (\ref{eq:Smat}) and the divergence of stresses is given by
$$
\Dop=\begin{pmatrix}g_x & g_y&0\\ 0 &g_x & g_y\\
-g_x^* & -g_y^*&0\\ 0 &-g_x^* & -g_y^* \end{pmatrix}.
$$
It differs from its counterpart in section \ref{sec:p1} because here it acts from vertices to cells.
\begin{figure}[t]
\centering
  \includegraphics[width=12cm]{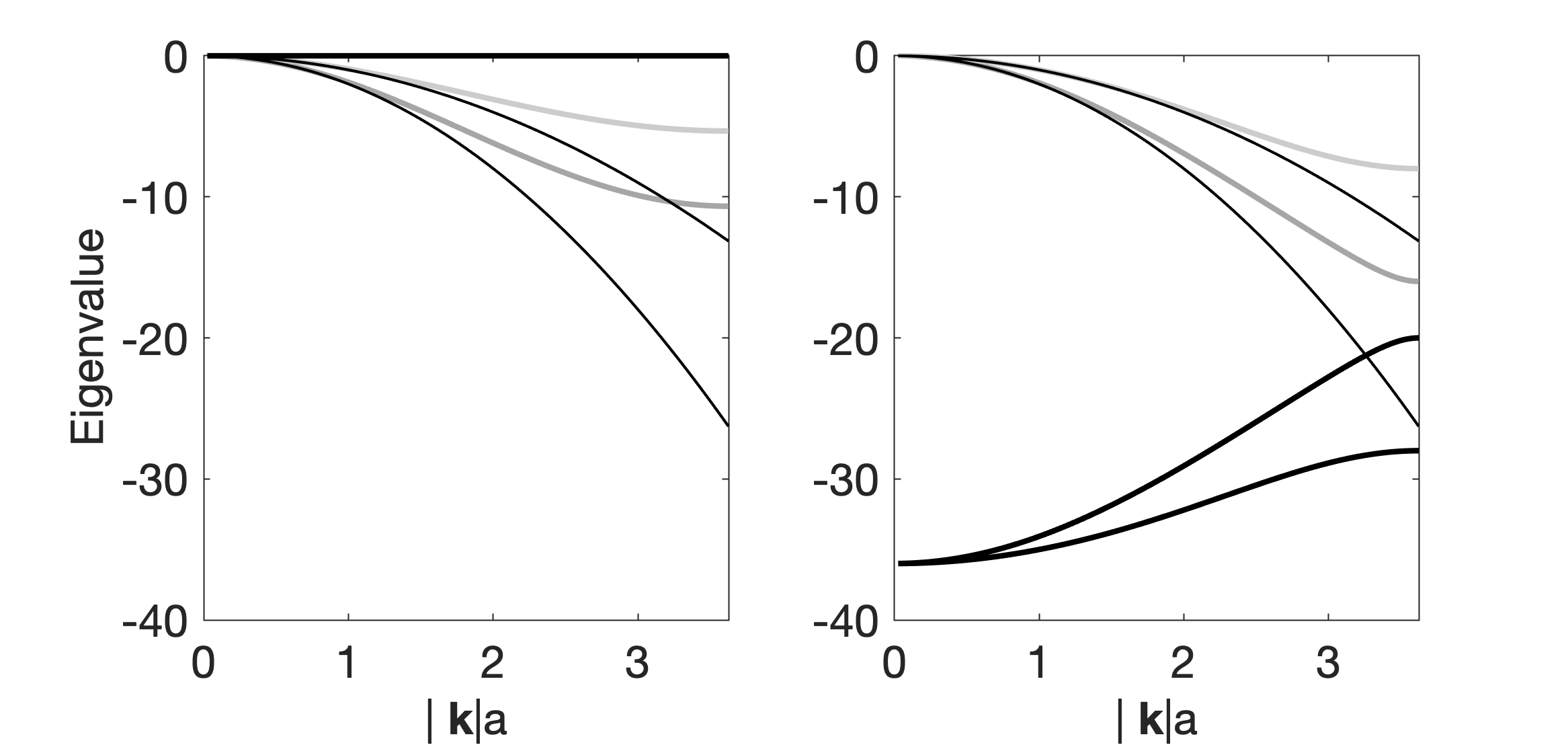}
\vspace{-0.5cm}\caption{The physical (thick gray lines) and spurious (thick black lines) eigenvalues of the Fourier symbol of $a^2\eta^{-1}\nabla\cdot\bsigma$  for the cell velocity placement as a function of $|\mathbf{k}|a$ for the wavevector oriented at $\pi/6$ to the $x$-axis and and $\zeta=\eta$. The thin black lines correspond to the dimensionless eigenvalues $-(k^2+l^2)a^2$ and $-2(k^2+l^2)a^2$ of the continuous case. Left panel: Strain rates are computed at vertices; Right panel: The proposed algorithm.}
\label{fig:p0}
\end{figure}

The system matrix is the product $\Dop\Sop\Eop$. The eigenvalues of $a^2\eta^{-1}\Dop\Sop\Eop$ are shown in the left panel of Fig. \ref{fig:p0}. It can be seen that the two physical eigenvalues (gray lines) are indeed reproduced if $|\mathbf{k}a|$ is small enough, and that they are slightly more accurate than in the case of vertex velocities (Fig. \ref{fig:p1}, left panel). However, in addition, there are two modes with zero eigenvalues (coinciding thick black lines). These modes are spurious. Their eigenvectors are spanning the null space of $\Dop\Sop\Eop$, which is the consequence of kernel in strain rates. They correspond to a grid-scale pattern of oscillations between $u$ and $d$ velocities and occur because $u$ and $d$ triangles are different. 

Case C is analyzed similarly. First, velocity derivatives are computed at triangles. For a $u$ triangle the components of velocity gradient are defined by three $d$ velocities on neighboring triangles, and vice versa. One readily sees that $-g_x^*$ and $-g_y^*$ will be involved in computations of strain rates on $u$ triangles, and $g_x$ and $g_y$  on $d$ triangles. The stresses are averaged to mid-edges and then divergence theorem is applied separately on $u$ and $d$ triangles to compute stress divergence. We omit further details here except for pointing out that case C shares the drawback of case V, and has less accurate physical branches, despite the fact that the dimension of the space where strain rates are first computed is twice as large as in case V.   

To eliminate the null space the strain rates and stresses need to be computed directly at edges such that no averaging is involved. This can be done in the least squares way. Alternatively, the strain rates of the case V (or C), averaged to edges, need to be 'corrected' to increase the weight of neighbor velocities. The following algorithm has been adopted. First, the velocity derivatives are computed at vertices and averaged to edges by taking, for every edge $e$, a half sum of values at its two vertices $v_1$ and $v_2$,
$$
(\partial_iu_j)^\star_e=(1/2)((\partial_iu_j)_{v_1}+(\partial_iu_j)_{v_2}).
$$

The 'corrected' edge velocity derivatives $(\partial_iu_j)_e$ are sought by minimizing
\begin{equation}
\label{eq:J}
J=((\partial_iu_j)^\star_e-(\partial_iu_j)_e)^2+\lambda(r_i(\partial_iu_j)_e-[u_j]_e).
\end{equation}
Here, $\lambda$ is the Lagrange multiplier, $\mathbf{r}$ is the vector connecting the centers of triangles sharing $e$ and $[\mathbf{u}]_e$ is the difference of velocities between these centers in the direction of $\mathbf{r}$. In Fig. \ref{fig:schematic}, for edge $e_1$,  $\mathbf{r}$ is the vector drawn from $c_6$ to $c_1$ and $[\mathbf{u}]_{e_1}=\mathbf{u}_{c_1}-\mathbf{u}_{c_6}$. Solving this minimization problem gives
$$
(\partial_x u)_e=(\partial_x u)^\star_e-\frac{r_x}{|\mathbf{r}|^2}(r_x(\partial_x u)^\star_e+r_y(\partial_y u)^\star_e-[u]_e),
$$
$$
(\partial_y u)_e=(\partial_y u)^\star_e-\frac{r_y}{|\mathbf{r}|^2}(r_x(\partial_x u)^\star_e+r_y(\partial_y u)^\star_e-[u]_e),
$$
for the derivatives of $u$ and similarly for the derivatives of $v$. The result reveals the motivation behind the correction: we replace the projection of velocity gradient on the vector $\mathbf{r}$ connecting velocity points across edge $e$ with its local estimate, i.e., the difference of velocity divided by the distance, $[\mathbf{u}]_e/|\mathbf{r}|$. 

After the derivatives are corrected and strain rates and stresses are computed at mid-edges, the stress divergence is computed in the finite-volume sense on triangles. The result is shown in the right panel of Fig. \ref{fig:p0}. A dramatic difference between the right and the left panels of Fig. \ref{fig:p0} is that the modes with zero eigenvalue in the left panel are replaced by the modes whose eigenvalues do not tend to zero at $|\mathbf{k}|a\to 0$ (black lines). These modes are still spurious because they do not tend to zero in the limit of small wavenumbers. In contrast to cases C and V, these modes are not expected to distort sea-ice dynamics because they have anomalously high negative eigenvalues even if $|\mathbf{k}|a\to 0$ and will decay fast. However, this has implications for the stability of internal time stepping in explicit methods such as EVP or mEVP as discussed in section \ref{sec:ill} below.

Note that despite the presence of spurious modes the accuracy of the representation of physical modes is improved compared to the case of vertex velocities (compare Fig \ref{fig:p0} and Fig. \ref{fig:p1}).

\subsection{Edge placement of velocities}
As proposed in \cite{Mehlmann2021}, the finite element method is used to represent velocities,   
$$
\mathbf{u}=\sum_{e\in\mathcal{E}}N_e\mathbf{u}_e,
$$
where $N_e$ is the non-conforming linear (Crouzeix--Raviart) basis function that equals 1 at edge $e$ and $-1$ at the vertex opposing $e$, and $\mathcal{E}$ is the set of mesh edges.  
Since the basis functions are linear, strain rates are constant on triangles. We associate their locations with triangle centers. For a  $u$-triangle, the Fourier amplitude of velocity gradient $\nabla \mathbf{u}$ is
$$
(\overline{\nabla \mathbf{u}})_c=(1/h)[\alpha(0,-2)\overline{\mathbf{u}}^a+\beta(-\sqrt{3},1)\overline{\mathbf{u}}^b+\gamma(\sqrt{3},1)\overline{\mathbf{u}}^c],
$$
where the phase multipliers are $\alpha=e^{-ilh/3},\beta=e^{-ika/4+ilh/6}$ and $\gamma=e^{ika/4+ilh/6}$. The phase differences featuring in $\alpha,\beta,\gamma$ correspond to phase shifts between the mid-edges and center of a $u$-triangle. They are twice smaller than the phase shifts in $\gamma_1^*, \beta_1^*, \alpha_1^*$ of section \ref{sec:p1}, which already indicates that one might get a better approximation in this case compared to the cases of vertex and cell velocities. 

The matrix connecting the vector of strain rate amplitudes $(\overline{\dot\epsilon_{xx}}^u,\overline{\dot\epsilon_{xy}}^u, \overline{\dot\epsilon_{yy}}^u, \overline{\dot\epsilon_{xx}}^d, \overline{\dot\epsilon_{xy}}^d, \overline{\dot\epsilon_{yy}}^d)^T$
with the vector of velocity amplitudes $(\overline{u}^a,\overline{v}^a,\overline{u}^b,\overline{v}^b,\overline{u}^c, \overline{v}^c)^T$ is written as
$$
\Eop=\begin{pmatrix}\Eop^u\\-(\Eop^u)^*\end{pmatrix}; \quad \Eop^u=\frac{1}{h}\begin{pmatrix}0&0&-\sqrt{3}\beta&0&\sqrt{3}\gamma&0\\ -\alpha&0&\beta/2&-\sqrt{3}\beta/2&\gamma/2&\sqrt{3}\gamma/2\\ 0&-2\alpha&0&\beta&0&\gamma\end{pmatrix}.
$$
We will get Fourier amplitudes of stresses by multiplying the strain rate amplitudes with $\Zop$ from (\ref{eq:Smat}). 

The divergence of stresses is computed in a weak sense as (see also \cite{Mehlmann2021}) 
\begin{equation}
\int \mathbf{N}_e\cdot(\Vop\mathbf{u}) dS=-\int\nabla\mathbf{N}_e:\bsigma dS -\varepsilon(2\eta/\ell_e)\int[\mathbf{N}_e]\cdot[\mathbf{u}] dl.
\label{eq:divnc}
\end{equation}
Here, $\mathbf{N}_e=\mathbf{w}_eN_e$ is the non-conforming linear test function with the amplitude equaled to amplitude of $\mathbf{w}_e$, and $\ell_e$ is the length of edge $e$. The stress divergence $\Vop\mathbf{u}$ is expanded as 
$$
\Vop\mathbf{u}=\sum_{e'\in\mathcal{E}}N_{e'}(\Vop\mathbf{u})_{e'}.
$$
The mass matrix with entries $\int N_eN_{e'}dS$ appearing on the left hand side of (\ref{eq:divnc}) is diagonal. Its diagonal entry is the area $A_e$ associated to edge $e$. For an equilateral triangular mesh it equals 2/3 of triangle area, and is three times smaller than the area of median-dual control volume in the case of vertex velocities. The second term on the right hand side is the stabilization proposed in  \cite{Mehlmann2021}  with $\varepsilon$ the stabilization parameter. This form of stabilization has been introduced by \cite{HansboLarson} for a linear elastic problem. The quantities in the square brackets are jumps across edges in the direction normal to the edges. The result is independent of the normal vectors, but the same orientation should be used while computing both differences. 

We will consider the contributions from two terms on the right hand side of (\ref{eq:divnc}) separately. To compute the first one, we need the matrix of the divergence operator, which is written as
$$
\Dop=-\frac{3}{2h}\begin{pmatrix}0&-2\alpha^*&0&0&2\alpha&0\\ 0&0&-2\alpha^*&0&0&2\alpha\\
-\sqrt{3}\beta^*&\beta^*&0&\sqrt{3}\beta&-\beta&0\\
0& -\sqrt{3}\beta^*&\beta^*&0&\sqrt{3}\beta&-\beta\\
\sqrt{3}\gamma^*&\gamma^*&0&-\sqrt{3}\gamma&-\gamma&0\\
0& \sqrt{3}\gamma^*&\gamma^*&0&-\sqrt{3}\gamma&-\gamma\end{pmatrix}.  
$$
It already contains the division by $A_e$ so that the matrix $\Dop\Zop\Eop$ corresponds to the Fourier symbol for unstabilized $\Vop$. The stabilization will add the matrix  $\varepsilon\Top$,
$$
\Top=\frac{2\eta}{ah} \Lop ^2,\quad 
\Lop =\begin{pmatrix}0&0&c_1&0&-c_2&0\\0&0&0&c_1&0&-c_2\\
		 -c_1&0&0&0&c_3&0\\ 0&-c_1&0&0&0&c_3\\
		 c_2&0&-c_3&0&0&0\\0&c_2&0&-c_3&0&0\end{pmatrix},
$$
where $c_1=2\cos(-ka/4+lh/2)$, $c_2=2\cos(ka/4+lh/2)$ and $c_3=2\cos(ka/2)$. Although the stabilization term $\Top$ has a structure of a discrete Laplacian operator, it has no 'physical' eigenvalues. Two of its branches have zero eigenvalues, and the others have finite eigenvalues even in the limit of vanishing wavenumbers, implying that the stabilization term works to couple edge velocities at grid scales. 

\begin{figure}[t]
\centering
  \includegraphics[width=12cm]{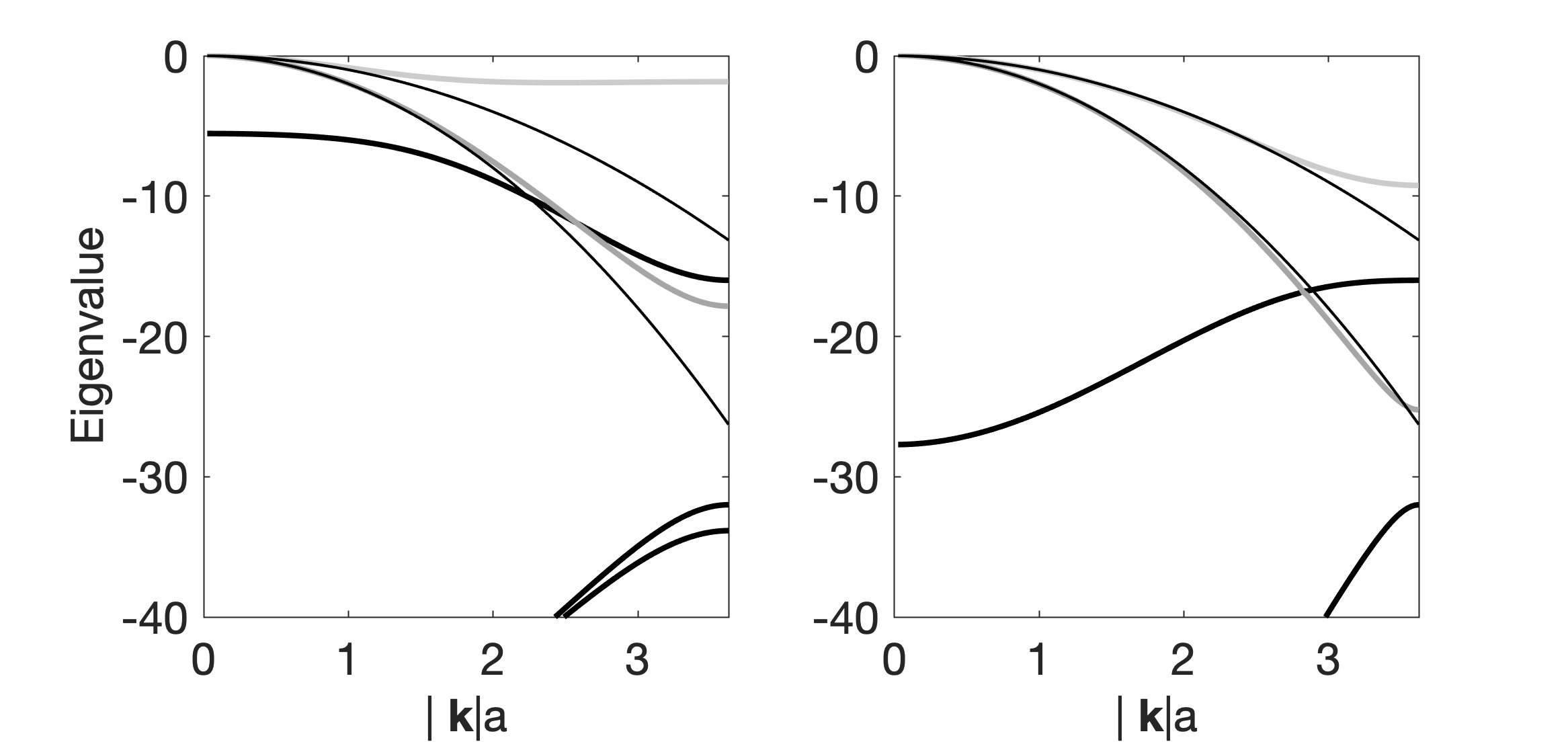}
\vspace{-0.5cm}\caption{The physical (thick gray lines) and spurious (thick black lines) eigenvalues  of the Fourier symbol of $a^2\eta^{-1}\nabla\cdot\bsigma$ for edge velocity placement as a function of $|\mathbf{k}|a$ for $\varepsilon=0.2$ (left panel) and $\varepsilon=1.0$ (right panel). The wavevector is at $\pi/6$ to the $x$-axis, and $\zeta=\eta$. The thin black lines plot $-(k^2+l^2)a^2$ and $-2(k^2+l^2)a^2$.}
\label{fig:pnc}
\end{figure}

We begin with noting that while the operators $\eta\Delta \mathbf{u}$ and $\nabla\cdot\bm{\sigma}$ coincide in the continuous case for $\zeta=0$, this does not happen in the case of non-conforming linear functions. If $\varepsilon=0$, the eigenvalues of $\nabla\cdot\bm{\sigma}$ (not shown) contain zero branches reflecting the presence of a non-trivial kernel in $\Eop$. However, the remaining branches do not approximate $-\eta(k^2+l^2)$, whereas the discrete $\eta\Delta \mathbf{u}$ approximates them. The lack of approximation persists for $\zeta\ne 0$. This behavior is related to a too compact stencil of the divergence operator: the first term in (\ref{eq:divnc}) accounts only for the derivative in the direction normal to edge for elementwise-constant stresses.  

The situation changes when stabilization is added, and the central point here is the magnitude of the stabilization coefficient. If the coefficient $\varepsilon$ is small, the physical branches of $\Vop$ start to emerge, as is seen in the left panel of Fig. \ref{fig:pnc} for $\varepsilon=0.2$ (the gray curves). However, they are accurate only for the small part of the range of resolved wavenumbers. The thick black curves correspond to spurious branches. They have similar origin as in the case of cell velocities: Only subsets of edges are invariant with respect to mesh translations, so there are modes that correspond to grid-scale oscillations. They will decay fast in viscous regime of sea ice dynamics because their eigenvalues are large in magnitude. Similarly to the case of cell velocities, there is an implication for the stability of time stepping of EVP-like methods (see below). 

The accuracy increases substantially if $\varepsilon$ is increased, as can be concluded from the right panel of Fig. \ref{fig:pnc}, which correspond to $\varepsilon=1$. There is some sensitivity to the orientation of $\mathbf{k}$, however up to $\varepsilon$ about 0.5 the gray curves approach asymptotically the theoretical curves (thin black lines) over an increasing range of wavenumbers if $\varepsilon$ increases. Further increase of $\varepsilon$ extends the range of wavenumbers where the curves are close to the theoretical ones, however the accuracy might be limited, as seen in Fig \ref{fig:pnc} for the lower gray branch for $|\mathbf{k}|a$ exceeding 2.      
In the end, the operator $\Dop\Zop\Eop+\varepsilon \Top$ is stable and its physical branches provide a good approximation of the Fourier symbol of continuous operator. Comparing the right panels of Fig. \ref{fig:p0} and \ref{fig:pnc} we can conclude that the eigenvalues can be more accurate than for the cell velocities if the stabilization parameter is properly selected. Counterintuitively, the accuracy of physical branches of $\Dop\Zop\Eop+\varepsilon \Top$ increases if $\varepsilon$ is increased to $\varepsilon>0.5$ or even higher. The rather accurate representation of physical branches does not imply that one can rely on the results of simulations on grid scales because spurious branches intersect with physical ones at $|\mathbf{k}|a$ between 2 and 3; they will contaminate solutions.

\section{Test case illustrations}
\label{sec:ill}
We use the test case described in \cite{Mehlmannetal2021} to illustrate the behavior of vertex-, cell- and edge-placed velocities and to show that their resolving capacity with respect to LKFs largely follows the insight derived from the analysis above. The test case explores the reaction of thin ice in a rectangular box to a cyclone moving along a box diagonal. Ice breaks and multiple LKFs are formed during the 2 days of simulations. The simulations are performed with  the mEVP method (\cite{Bouillon2013}) on a mesh made of equilateral triangles with the side of 2 km covering a square area 512 by 512 km in size. The sea-ice component of FESOM is used, extended for triangular B- and CD-grid discretizations. 
\begin{table}[t]
    \begin{center}
    \begin{tabular}{c|c|c|c}
    velocity placement& A-grid& B-grid & CD-grid  \\
  \hline
    $\alpha=\beta$&  500& 1200 & 1500\\ 
    \end{tabular}
    \caption{  Minimum values of the mEVP stabilization parameters $\alpha, \beta$ to achieve stable simulations}
    \label{stabi}
     \end{center}
\end{table}

We begin with a brief discussion of the time step limitation. According to \cite{Kimmritz2015}, for numerical stability the parameters $\alpha$ and $\beta$ of mEVP should be selected such that 
\begin{align}\label{Kimmritz:stabi}
  \alpha\beta >>
\Lambda^2\zeta\Delta t/m,  
\end{align}
where $-\Lambda^2$ is the largest negative eigenvalue of the Fourier symbol of $\eta^{-1}\nabla\cdot\bsigma$ operator and $m$ is the ice mass per unit area. As can be concluded by comparing the left panel of Fig. \ref{fig:p1} to the right panel of Fig. \ref{fig:p0}, the maximum of $\Lambda^2$ for cell velocities exceeds that for vertex velocities by a factor around 3.5. The maximum for edge velocities (not shown in Fig. \ref{fig:pnc}) is about 7 times higher. Although all conclusions of linear analysis above, as well as the stability analysis in \cite{Kimmritz2015}, cannot be accurate in a full nonlinear case, they remain qualitatively correct, as illustrated below. The larger $\Lambda^2$ implies that either the external time step $\Delta t$ must be reduced, or $\alpha$ and $\beta$ must be increased for stability if cell and edge velocities are used compared to the case of vertex velocities. 
To illustrate qualitative agreement with theory we performed simulations on A- and CD-grids selecting $\alpha=\beta=1000$ and $N_{EVP}=100$, and varying the external time step $\Delta t$. For A-grid, runs with $\Delta t=2, 4, 7.5$ min showed no traces of noise, but weak noise started to appear for $\Delta t=15$ min. For the CD-grid run with $\Delta t=3$ min showed strong noise. The noise became weak for $\Delta t=2$ min and almost (but not fully) disappeared for $\Delta t=1.5$ min. Relating the cases with weak noise or just no noise we see that the ratio of external time steps is about 5--7, which agrees with \cite{Kimmritz2015} and the eigenvalue increase predicted by the Fourier analysis. 

In practice the external time step will be defined by an ocean model, and the time step about $\Delta t$=2 min will be expected on a mesh with 2 km cells. In a line with (\ref{Kimmritz:stabi}) it has been found that the simulation is stable for the edge velocities if $\alpha=\beta=1500$, whereas they could be reduced to 1200 and 500 respectively for cell and vertex velocities. The choice of $\alpha, \beta$ is collected in Table \ref{stabi}. For the CD-grid we tested that the sensitivity of the LKF pattern to the magnitude of $\alpha,\beta$ in the range (1500-2500) and the number of substeps $N_{EVP}$ in the range 100-2000 is rather weak compared to the differences created by using different locations for velocities. Similar weak sensitivity is observed for other discretizations. Larger $\alpha, \beta$ imply that mEVP solutions may further diverge from VP solutions unless $N_{EVP}$ is high enough.     

\begin{figure}[t]
\centering
  \includegraphics[width=12cm]{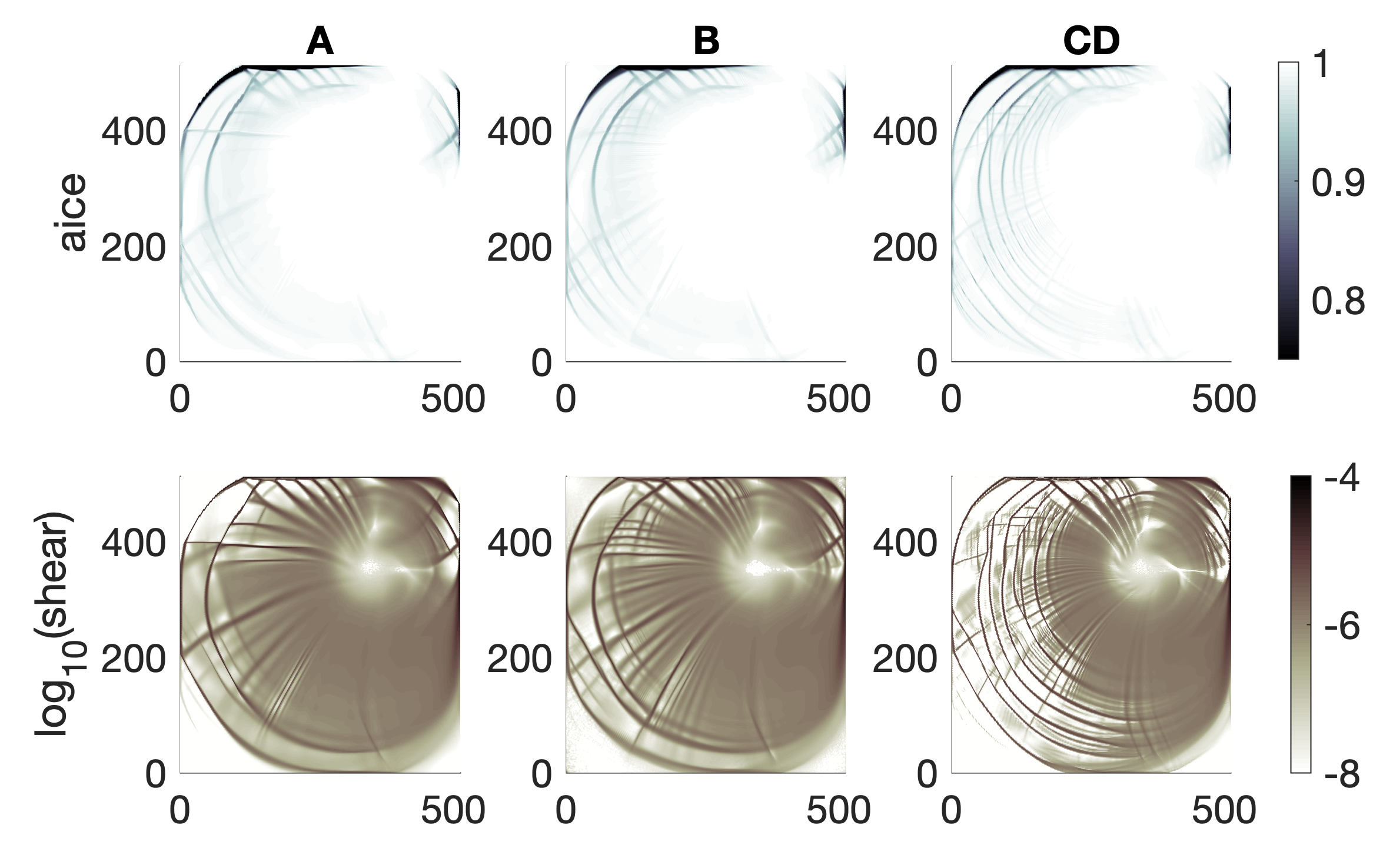}
\vspace{-0.5cm}\caption{The patterns of sea-ice concentration (top row) and shear (bottom row) for vertex (left column, A-grid), cell (middle column, B-grid) and edge (right column, CD-grid) velocity placement. The number of simulated LKFs increases together with the number of degrees of freedom in the velocity field. }
\label{fig:3as}
\end{figure}

Figure \ref{fig:3as} presents the patterns of sea ice concentration and shear defined as $((\varepsilon_{11}-\varepsilon_{22})^2+4\varepsilon_{12}^2)^{1/2}$, simulated in the test case of \cite{Mehlmannetal2021} with A-, B- and CD-grid discretizations using $\alpha,\beta$ defined in Table \ref{stabi}. On qualitative level, it is seen that the number of simulated LKFs increases together with an increase in the number of degrees of freedom in the velocity field, which also correlates with improved accuracy of representing the eigenvalues of linear operator of stress divergence by cell and edge discretizations.

The patterns of sea ice concentration are commonly 
smooth, and the performance of particular discretization in the context considered here should be judged by the behavior of strain rates which amplify grid-scale noise. The pattern of shear for vertex velocities (bottom left panel) is clean, as expected. But such patterns contain no obvious noise also for cell and edge discretizations. This indicates that the stabilization proposed in \cite{Mehlmann2021} for edge velocities and the one proposed here for cell velocities are efficient in removing the kernels in strain rates or stress divergence and that spurious modes do not create problems other than time step limitation. 

Some minor issues can still be seen on closer inspection. The middle bottom panel shows weak small-amplitude grid-scale noise in box western corners where sea ice velocities are nearly zero. This noise originates from boundaries and is by all probability related to the fact that no-slip boundary conditions are implemented only approximately because all velocity points are inside the box. It  weakens if replacement pressure is not used, which is the case in the middle column of Fig. \ref{fig:3as}. We hope that taking boundary at velocity points will eliminate this behavior.  

For edge velocities, some of the LKFs seen in the strain field are only one grid cell apart, which are not the scales one may rely on. This signals that further fine-tuning of the amplitude of stabilization might be needed, which was not attempted here.

\section{Discussions}

Cell (B-grid) and edge (CD-grid) placement of sea ice velocities results in a higher number of simulated LKFs compared to the vertex (A-grid) placement. This correlates with a larger number of DoF on B- and CD-grids as there are more triangles and edges than vertices on triangular meshes. This also correlates with a higher accuracy of the discrete operators on B- and CD-grids as demonstrated in section \ref{sec:fourier}. However, neither cell nor edge sea ice velocity placement leads to a straightforward discretization. Both need special measures to avoid a kernel in the discrete strain rate or stress divergence operators. The kernels are eliminated if appropriate measures are taken. The Fourier analysis, despite its obvious limitations, demonstrates the consequences of these measures. 

For cell velocities, the strategy of eliminating a kernel in stress divergence becomes clear after learning that the kernel is created either through rank deficiency (case V, vertex locations of strain rates) or averaging (case C, cell location of strain rates). The remaining choice is between direct computation or correction of strain rates on edges. We used a correction. Its key element is the increase in the contribution from the differences between across-edge velocities. The idea is analogous to an implicit discretization of gradients suggested in expression (18) of \cite{Hutchings2004}. The Fourier analysis directly shows that this not only eliminates the kernel, but also improves accuracy. The strong constraint in the minimization problem (\ref{eq:J}), introduced though the Lagrangian multiplier $\lambda$, can be relaxed to a weak constraint $(w/|\mathbf{r}|)^2(r_i(\partial_iu_j)_e-[u_j]_e)^2$, where $w$ is the dimensionless weight. We, however, found that only rather large weights ensure noise-free behavior, so that the strong constraint in (\ref{eq:J}) seems to be a more reliable option.   

For edge velocity an important insight from the Fourier analysis is that the stabilization ensures approximation: without it the first term in (\ref{eq:divnc}) not only has a kernel, but is also not able to approximate the eigenvalues of continuous problem. Common appreciation is that stabilization has to be kept as weak as possible in order not to affect physical branches. Here it has to be kept sufficiently strong in order to ensure that the physical branches are recovered. In the EVP or mEVP methods, the stresses $\bsigma$ in (\ref{eq:divnc}) differ from those of the VP method, because they are computed in a time stepping or iterative procedure. This raises a new question on the suitable amplitude of $\varepsilon$ (see \cite{Mehlmann2021} for a practical choice). Furthermore, the stabilization term in (\ref{eq:divnc}) becomes an additional factor influencing numerical stability of EVP-like methods and taking $\varepsilon$ that is too large may affect it. The insight that it should be sufficiently large remains valid, but the selection of optimal $\varepsilon$ for EVP and mEVP methods requires additional studies. 

Both cell and edge velocity placements maintain spurious modes. These modes acquire anomalously large negative eigenvalues when kernels in operators are eliminated, as illustrated in Fig. \ref{fig:p0} and partly in Fig. \ref{fig:pnc}. Their origin is rooted in the geometry of triangular meshes, as has been briefly explained above (see \cite{DanilovKutsenko2019} for more detail). Because one half or two thirds of modes are spurious, respectively, for discretizations based on cell and edge velocities, there might be an impression that the velocity DoF are used suboptimally in these cases. However, despite spurious modes, these discretizations ensure a much improved accuracy of the representation of physical modes compared to the A-grid physical modes, which correlates with their substantially increased resolving capability. An interesting question for future studies is how the A-grid discretization would compare to B- and CD-grid discretizations on meshes with a matching number of DoF.  

The only apparent consequence of spurious modes seen in our test case simulations is their impact on the stability of explicit EVP-like methods. The stability conditions are tighter for cell and edge velocities because spurious modes are characterized by anomalously high negative eigenvalues. We hope that this will be the only consequence in realistic configurations, but it remains to be seen.

\section{Conclusions} 

Elimination of kernels in discrete stress divergence is a key requirement for numerical stability of discretizations of sea ice dynamics on triangular meshes using cell- or edge-based velocities. This is achieved through the stabilization proposed by \cite{Mehlmann2021} for edge velocities and the procedure proposed here for cell velocities. 

In both cases the consequence of stabilization is that dimensionless eigenvalues of spurious modes, supported on B- and CD-grids, take large negative values. As a result, spurious modes are not expected to distort solutions: if excited, they will decay faster than physical branches. However, this will impact stability of time stepping in explicit methods, as discussed in section \ref{sec:ill}. It should be expected that B- and CD-grids will require smaller internal time steps in EVP or larger stability parameters in mEVP.     

Summing up, the cell- and edge-based sea ice velocities seem to be a promising alternative to the vertex placement for discretizing equations of sea ice dynamics on triangular meshes. They ensure an improved accuracy in representing physical modes of stress divergence operator and higher effective resolution. It remains to be seen how well cell- and edge-based discretizations perform in realistic conditions, how sea ice dynamics are affected by the placement of scalar degrees of freedom and what is an optimal stabilization for the edge velocities in explicit time integration methods. 

It is hoped that insights provided by the simple approach in this work will be helpful for modelers working with sea ice dynamics on unstructured triangular (or their dual) meshes.

\section*{Acknowledgments} This work is a contribution to project S2 of the Collaborative Research Centre TRR181 "Energy Transfer in Atmosphere and Ocean" funded by the Deutsche Forschungsgemeinschaft (DFG, German Research Foundation) - Projektnummer 274762653.


\begin{thebibliography}{10}

\bibitem{Arakawa1977}
A.~Arakawa and V.R. Lamb.
\newblock Computational design of the basic dynamical processes of the {UCLA}
  general circulation model.
\newblock {\em Methods Comput. Phys.}, 17:173–265, 1977.

\bibitem{Bouillon2013}
S.~Bouillon, T.~Fichefet, V.~Legat, and G.~Madec.
\newblock {T}he elastic-viscous-plastic method revisited.
\newblock {\em Ocean Modelling}, 71:2--12, 2013.

\bibitem{Coon2007}
Max Coon, Ron Kwok, Gad Levy, Matthew Pruis, Howard Schreyer, and Deborah
  Sulsky.
\newblock Arctic ice dynamics joint experiment {(AIDJEX)} assumptions revisited
  and found inadequate.
\newblock {\em Journal of Geophysical Research: Oceans}, 112(C11), 2007.

\bibitem{DanilovKutsenko2019}
S.~Danilov and A.~Kutsenko.
\newblock On the geometric origin of spurious waves in finite-volume
  discretizations of shallow water equations on triangular meshes.
\newblock {\em J. Comput. Phys.}, 398:108891, 2019.

\bibitem{Danilov2017}
Sergey Danilov, Dmitry Sidorenko, Qiang Wang, and Thomas Jung.
\newblock The {F}inite-volum{E} {S}ea ice--{O}cean {M}odel ({FESOM2}).
\newblock {\em Geosci. Model Dev.}, 10:765--789, 2017.

\bibitem{Danilov2015}
Sergey Danilov, Qiang Wang, Ralph Timmermann, Nikolay Iakovlev, Dmitry
  Sidorenko, Madlen Kimmritz, Thomas Jung, and Jens Schr{\"o}ter.
\newblock Finite-element sea ice model ({FESIM}), version 2.
\newblock {\em Geoscientific Model Development}, 8(6):1747--1761, 2015.

\bibitem{Feltham08}
D.L. Feltham.
\newblock Sea {I}ce {R}heology.
\newblock {\em Annual Review of Fluid Mechanics}, 40:91--112, 2008.

\bibitem{Gao2011}
G.~Gao, C.~Chen, J.~Qi, and R.~C. Beardsley.
\newblock An unstructured-grid, finite-volume sea ice model: {D}evelopment,
  validation, and application.
\newblock {\em J. Geophys. Res.}, 116:C00D04, 2011.

\bibitem{HansboLarson}
P.~Hansbo and M.~G. Larson.
\newblock Discontinuous {G}alerkin and the {C}rouzeix-{R}aviart element:
  {A}pplication to elasticity.
\newblock {\em Mathematical Modelling and Numerical Analysis}, 37:63--72, 2003.

\bibitem{Hibler1979}
W.~D. Hibler, III.
\newblock A {D}ynamic {T}hermodynamic {S}ea {I}ce {M}odel.
\newblock {\em J. Phys. Oceanogr.}, 9:815--846, 1979.

\bibitem{HunkeDukowicz1997}
Elizabeth~C. {Hunke} and J.~K. Dukowicz.
\newblock An {E}lastic-{V}iscous-{P}lastic model for sea ice dynamics.
\newblock {\em J. Phys. Oceanogr.}, 27:1849--1867, 1997.

\bibitem{Hutchings2004}
J.~K. Hutchings, H.~Jasak, and S.~W. Laxon.
\newblock A strength implicit correction scheme for the viscous-plastic sea ice
  model.
\newblock {\em Ocean Modelling}, 7:111--133, 2004.

\bibitem{Hutter2018a}
N.~Hutter, L.~Zampieri, and M.~Losch.
\newblock Leads and ridges in arctic sea ice from rgps data and a new tracking
  algorithm.
\newblock {\em The Cryosphere}, 13(2):627--645, 2019.

\bibitem{Kimmritz2015}
M.~Kimmritz, S.~Danilov, and M.~Losch.
\newblock On the convergence of the modified elastic-viscous-plastic method for
  solving the sea ice momentum equation.
\newblock {\em J. Comp. Phys.}, 296:90--100, 2015.

\bibitem{Korn2017}
P.~Korn.
\newblock Formulation of an unstructured grid model for global ocean dynamics.
\newblock {\em J. Comput. Phys.}, 339:525--552, 2017.

\bibitem{Lietaer2008}
O.~Lietaer, T.~Fichefet, and V.~Legat.
\newblock The effects of resolving the {C}anadian {A}rctic {A}rchipelago in a
  finite element sea ice model.
\newblock {\em Ocean Modelling}, 24:140--152, 2008.

\bibitem{Mehlmannetal2021}
C.~Mehlmann, S.~Danilov, M.~Losch, J.-F. Lemieux, N.~Hutter, T.~Richter,
  P.~Blain, E.~C. Hunke, and P.~Korn.
\newblock Simulating linear kinematic features in viscous-plastic sea ice
  models on quadrilateral and triangular grids.
\newblock {\em http://arxiv.org/abs/2103.04431}, 2021.

\bibitem{Mehlmann2021}
C.~Mehlmann and P.~Korn.
\newblock Sea-ice dynamics on triangular grids.
\newblock {\em J. Comput. Phys.}, 428:110086, 2021.

\bibitem{Petersen2019}
M.~R.and Asay-Davis Petersen, X.~S., Q.~Berres, A. S.and~Chen, N.~Feige, M.~J.
  Hoffman, D.~W. Jacobsen, P.~W. Jones, M.~E. Maltrud, S.~F. Price, T.~D.
  Ringler, G.~J. Streletz, A.~K. Turner, L.~P. Van~Roekel, M.~Veneziani, J.~D.
  Wolfe, P.~J. Wolfram, and J.~L. Woodring.
\newblock An evaluation of the ocean and sea ice climate of {E3SM} using {MPAS}
  and interannual {CORE-II} forcing.
\newblock {\em J. of Advances in Modeling Earth Systems}, 11:1438–1458, 2019.

\bibitem{Ringler2013}
T.~Ringler, M.~Petersen, R.~Higdon, D.~Jacobsen, M.~Maltrud, and P.W. Jones.
\newblock A multi-resolution approach to global ocean modelling.
\newblock {\em Ocean Modell.}, 69:211--232, 2013.

\bibitem{Timmermann2009}
R.~Timmermann, S.~Danilov, J.~Schr\"oter, C.~B\"oning, D.~Sidorenko, and
  K.~Rollenhagen.
\newblock Ocean circulation and sea ice distribution in a finite element global
  sea ice - ocean model.
\newblock {\em Ocean Modell.}, 27, 2009.

\bibitem{Wang2016}
Q~Wang, S~Danilov, T~Jung, Lars Kaleschke, and A~Wernecke.
\newblock Sea ice leads in the {A}rctic {O}cean: {M}odel assessment,
  interannual variability and trends.
\newblock {\em Geophysical Research Letters}, 43(13):7019--7027, 2016.

\bibitem{Wang2014}
Qiang Wang, Sergey Danilov, Dmitry Sidorenko, Ralph Timmermann, Claudia
  Wekerle, Xuezhu Wang, Thomas Jung, and Jens Schr{\"o}ter.
\newblock The {F}inite {E}lement {S}ea {I}ce-{O}cean model ({FESOM}) v. 1.4:
  formulation of an ocean general circulation model.
\newblock {\em Geoscientific Model Development}, 7(2):663--693, 2014.

\end{thebibliography}
\end{document}